# A new formula for the n-th prime:

Sebastián Martín Ruiz
Avda. de Regla 43, Chipiona 11550 Spain
smruiz@telefonica.net

Using the below expression for the Characteristic Function of Prime Numbers:

$$\left\lfloor \frac{lcm(1,2,...,j)}{j \cdot lcm(1,2,...,j-1)} \right\rfloor = \begin{cases} 1 & \text{if } j \text{ is prime} \\ 0 & \text{if } j \text{ is composite} \end{cases}$$

(This function is complementary to the Smarandache Prime Function [1], defined as:
P(n) = 0, if n is prime, and P(n) = 1 otherwise.
It is easy to prove this expression studying it in detail.
 I have obtained this expression this last month (March 2004) but I do not know if already it is known.)

we obtain the following formula for the n-th prime [2],[3],[4],[5],[6]:

$$p_n = 1 + \sum_{k=1}^{\lfloor 2n \log n + 2 \rfloor} \left( 1 - \left\lfloor \frac{1}{n} \sum_{j=2}^{k} \left\lfloor \frac{lcm(1,2,...,j)}{j \cdot lcm(1,2,...,j-1)} \right\rfloor \right\rfloor \right)$$

The Proof is the same of the previous articles.
It is necessary to see the references for a complete comprehension of the formulas.

We can see that this formula is faster than the previous:

Comparative table of times:

| Prime    | Prob 38   | Lcm       | Prob 38 mod | Lcm mod   |
|----------|-----------|-----------|-------------|-----------|
| P10=29   | 0,2   sec | 0,06  sec | 0,04 sec    | 0,02  sec |
| P20=71   | 2,5   sec | 0,8   sec | 0,4   sec   | 0,1   sec |
| P30=113  | 12,2  sec | 3,7   sec | 1,3   sec   | 0,4   sec |
| P40=173  | 36,5  sec | 10,5  sec | 2,9   sec   | 0,7   sec |
| P50=229  | 84    sec | 24    sec | 5,5   sec   | 1,3   sec |
| P100=541 |           |           | 41    sec   | 6,9   sec |
| P200=1223|           |           | 299   sec   | 39    sec |

*Prob 38: [6] It is the original formula without modifying.
The time complexity of this algorithm is  O(nlog n)^3
*Lcm: Is the new formula.
*Prob 38 Mod: [6] It is the original formula with the modifications.
The time complexity of this algorithm is  O(n log n)^(3/2).
*Lcm mod: It is the new formula calculating lcm (1,2, …, j) of way recurrent.

**Which is the time complexity of this algorithm?**

The code in Mathematica:
```
L[1]=1;
L[n_]:=L[n]=LCM[L[n-1],n]
LG[n_]:=L[n]/L[n-1]
FL[n_]:=Quotient[LG[n],n]
Pii[n_]:=Sum[FL[i],{i,2,n}]
PrimeLCM[n_]:=1+Sum[1-Quotient[Pii[k],n],{k,1,Floor[2*n*Log[n]+2]}]
Do[Print[n," ", Timing[PrimeLCM[n]]," ",Prime[n]],{n,200,200}]
200  1223    1223
{39.438 Second}
```

We can accelerate it more enough of the following form:

Using the bound of Rosser and Schoenfeld for $p_n$ [7]:

$$c_n = n \log n + n(\log(\log n) - 1/2)$$

and modifying the formula considering that $p_n > \lfloor n \log n \rfloor$ [7] we obtain for n> 1:

$$p_n = \lfloor n \log n \rfloor + \sum_{k=\lfloor n \log n \rfloor}^{\lfloor C_n+3 \rfloor} \left(1 - \left\lfloor \frac{1}{n} \sum_{j=2}^{k} \left\lfloor \frac{lcm(1,2,...,j)}{j \cdot lcm(1,2,...,j-1)} \right\rfloor \right\rfloor \right)$$

This is a nice expression that relates the n-th prime number with the approximation obtained with the prime number theorem $n \log n$ adding a term of error.

The new times are:

| Prime | Prob 38 mod | Lcm mod | RS acceleration |
|---|---|---|---|
| P10=29 | 0,04 sec | 0,02 sec | 0 sec |
| P20=71 | 0,4 sec | 0,1 sec | 0,02 sec |
| P30=113 | 1,3 sec | 0,4 sec | 0,05 sec |
| P40=173 | 2,9 sec | 0,7 sec | 0,09 sec |
| P50=229 | 5,5 sec | 1,3 sec | 0,15 sec |
| P100=541 | 41 sec | 6,9 sec | 0,86 sec |
| P200=1223 | 299 sec | 39 sec | 4,59 sec |

**Which is the time complexity of this algorithm?**